\documentclass[12pt]{smfart}
\usepackage[frenchb]{babel}
\usepackage{mathptmx}
\usepackage{amsmath}
\usepackage{amssymb}
\usepackage{smfthm}
\usepackage{amscd}
\usepackage{smfenum}

\usepackage{amssymb,amsfonts}

\usepackage{diagrams}
\usepackage{graphics}
\usepackage[T1]{fontenc}     %fontes a rajouter pour accents
\usepackage[latin1]{inputenc}

\newcounter{ctrmasection}
\newcounter{ctrmasoussection}

\newenvironment{pf*}{  \bf{D{\'e}monstration : }\rm}{$\Box$}

{\begin{displaymath} \begin{array}{c}}%
{\end{array} \end{displaymath}}
\newcommand{\Aut}{\operatorname{Aut}}

\newcommand{\Pic}{\operatorname{Pic}}

\newcommand{\Card}{\operatorname{Card}}
\newcommand{\Poly}{\operatorname{Poly}}

\begin{document}
\title[Invariants g\'eom\'etriques et  automorphismes du plan affine ]{Sur des invariants g\'eom\'etriques associ\'es aux automorphismes du plan affine} 
\alttitle{Geometric invariants of automorphisms of the affine plane}
\author{Sandra Marcello}
%\address{Adresse permanente, Th{\'e}orie des Nombres, Institut de Math{\'e}matiques de Jussieu
%175, rue du Chevaleret, 75013 Paris, France\\}
%\email{marcello@math.jussieu.fr\\}

\address{Max-Planck-Institut f{\"u}r Mathematik,
Vivatsgasse 7, 53111, Bonn, Deutschland \\}
\email{marcello@mpim-bonn.mpg.de}
\date{\today}

\begin{abstract}Nous  associons \`a chaque 
 automorphisme du plan affine, une construction g\'eom\'etrique poss\'edant certaines propri\'et\'es, la {\it{r\'esolution canonique}}. Nous \'etudions la g\'eom\'etrie de la r\'esolution canonique, en d\'eduisons une majoration d'un invariant g\'eom\'etrique (l'indice ample) associ\'e \`a un automorphisme du plan affine et relions la structure du cone effectif de la surface obtenue aux valeurs d'un autre invariant g\'eom\'etrique (l'indice effectif).
\end{abstract}

\begin{altabstract} We associate to each automorphism of the plane, a geometric construction with some properties, it is the {\it{canonical resolution}}.
 We study the geometry of the canonical resolution, we deduce from it an upper bound for a geometric invariant (the ample index) linked to an automorphism of the affine plane and we link the structure of the effective cone of the surface we get to values of another geometric invariant (the effective index).
\end{altabstract}


\altkeywords{Automorphisms of affine plane,\  ample cone, \ effectice cone}

\maketitle

\tableofcontents

\mainmatter

\section{Introduction}

En 1994, dans le but d'obtenir des propri\'et\'es arithm\'etiques  de l'application de H\'enon quadratique, J. Silverman \cite {Sil1} \'etudie la g\'eom\'etrie de cette application. Dans ce texte nous d\'efinissons, gr\^ace \`a des propri\'et\'es intrins\`eques des automorphismes du plan affine, cette construction de mani\`ere g\'en\'erale pour tout automorphisme du plan affine et \'etudions les propri\'et\'es g\'eom\'etriques   de ce que nous appelons une {\it{r\'esolution canonique}}. La r\'esolution canonique nous permet de calculer ou de majorer des invariants g\'eom\'etriques (les indices amples  et effectifs) que nous avons d\'efinis dans \cite{moi4} dans le but d'obtenir des r\'esultats de nature arithm\'etiques. Pour plus de pr\'ecisions sur les indices amples et effectifs nous renvoyons le lecteur \`a  \cite{moi4} et pour une motivation des questions arithm\'etiques auxquelles nous nous int\'eressons, nous renvoyons le lecteur \`a \cite{1art}.
L'une de ces questions est un analogue aux probl\`emes classiques de d\'ecompte des  points rationnels sur les vari\'et\'es pour une synth\`ese voir \cite{Peyre}. 

Il est \`a noter qu'en vue d'applications arithm\'etiques il est int\'eressant de conna\^ \i tre la valeur de ces invariants g\'eom\'etriques.

Le corps de base est $\mathbb{C}$.
Nous notons $\mathbb{A}^2$ (resp. $\mathbb{P}^2$) l'espace affine (resp. projectif) de
dimension $2$ et
$H$ la droite {\`a} l'infini. Lorsque $\phi\in\Aut
(\mathbb{A}^2)$ est un automorphisme de $\mathbb{A}^2$, nous noterons encore
$\phi:\mathbb{P}^2\dots\rightarrow \mathbb{P}^2$ l'application rationnelle induite et d{\'e}signerons par
$Z(\phi)$ le {\it{lieu de non-d{\'e}finition}}, {\it{i.e}} le lieu g{\'e}om{\'e}trique, contenu dans $H$, o{\`u} $\phi$ n'est pas d{\'e}finie.

\noindent De plus, si $\phi\in\Aut(\mathbb{A}^2)$, alors $\phi$ est
d{\'e}fini par des polyn{\^o}mes $P_1,P_2$. On d{\'e}finit le {\it{degr{\'e}
    alg{\'e}brique}} ou  {\it{degr{\'e}}} de $\phi$ comme suit : $deg(\phi)=d=\max_i(d_i)$ avec $d_i$ degr{\'e} total de $P_i$.
Dans \cite{moi4}, nous avons construit les objets en dimension quelconque, nous nous contenterons ici de rappeller les d\'efinitions dans le cadre de la dimension 2.

\begin{defi}Soit V une vari{\'e}t{\'e} lisse projective. On dit que le morphisme birationnel $\pi:V\rightarrow\mathbb{P}^2$  est une {\it
r{\'e}solution} de $\phi$ (not{\'e}e  $(\pi,V)$) si $\phi\circ\pi$ et
$\phi^{-1}\circ\pi$ sont des morphismes de $V$ vers $\mathbb{P}^2$.
Parfois nous noterons:
$\psi=\phi\pi$ et $\psi'=\phi^{-1}\pi$.
\end{defi}

Dans ce texte nous d\'efinissons une r\'esolution particuli\`ere que nous appelerons {\it{r\'esolution canonique}}, gr\^ace \`a celle-ci nous majorons ou calculons des invariants g\'eom\'etriques associ\'es aux automorphismes du plan affine.
Ces invariants g\'eom\'etriques ont \'et\'e d\'efinis par l'auteure dans l'optique d'applications  arithm{\'e}tiques qui apparaissent clairement dans \cite{moi4} th{\'e}or{\`e}me A, 
ces invariants sont d\'efinis comme suit.

\paragraph*{Notations} 
\begin{itemize}
\item Si $V$ est une vari{\'e}t{\'e}, on note respectivement $\Pic_{\mathbb{Q}}(V)$ et $\Pic_{\mathbb{R}}(V)$ son
groupe de Picard tensoris{\'e} respectivement par $\mathbb{Q}$  et $\mathbb{R}$, $\Pic^+_{\mathbb{Q}}(V)$ et $\Pic^+_{\mathbb{R}}(V)$ les cone ferm{\'e}s
engendr{\'e}s par les classes de diviseurs effectifs, $\Pic^{a}_{\mathbb{R}}(V)$ le cone (ouvert)
engendr{\'e} par les classes de diviseurs amples et  $\Pic^{nef}_{\mathbb{R}}(V)$ le cone des diviseurs nef.
\item Soit $(\pi,V)$ une r\'esolution de $\phi$ et $\alpha\in\mathbb{R}$, nous notons:
$$D(\alpha,\pi,V):=(\phi\circ\pi)^*H+(\phi^{-1}\circ\pi)^*H-\alpha  \pi^*H.$$
\item Les nombres suivants sont \'egalement d\'efinis:
$$\alpha_{max,eff}(\phi,\pi,V):=\sup\{\alpha\in\mathbb{R}\;|\;D(\alpha,\pi,V) \in
\Pic^+_{\mathbb{R}}(V)\}, $$
$$\alpha_{max,amp}(\phi,\pi,V):=\sup\{\alpha\in\mathbb{R}\;|\;D(\alpha,\pi,V) \in
\Pic^{a}_{\mathbb{R}}(V)\}.$$
\end{itemize}

\begin{defi} Soit $\phi \in \Aut(\mathbb{A}^2)$. Les nombres suivants seront appel\'es respectivement l'{\it{indice effectif}} et l'{\it{indice ample}}:
$$\alpha(\phi,eff):=\sup\{\alpha_{\max, eff}(\phi,\pi,V)| \pi, V\},$$
$$\alpha(\phi,amp):=\sup\{\alpha_{\max,amp}(\phi,\pi,V)| \pi, V\}.$$
\end{defi}

Dans \cite{moi4}, nous montrons que ces nombres sont bien des invariants g{\'e}om{\'e}triques des automorphismes de l'espace affine.

\noindent Des propri{\'e}t{\'e}s propres aux automorphismes du plan affine nous permettent de d{\'e}finir la notion de {\it{r{\'e}solution canonique}} (voir sa construction paragraphe 4). Cette r{\'e}solution nous permet de d\'emontrer le th\'eor\`eme suivant:

\paragraph*{{\bf{Th\'eor\`eme A}}}
{\it{Soit $\phi\in \Aut(\mathbb{A}^2)$ de degr{\'e} alg{\'e}brique $d\geq
    2$. Alors, }}
$$\alpha(\phi, amp)\leq 0.$$

Dans \cite{moi4}, nous donnons une g\'en\'eralisation de ce th\'eor\`eme en dimension sup\'erieure avec toutefois une condition g\'eom\'etrique sur les automorphismes de l'espace affine.

Peu d'informations sont connues sur la structure du cone effectif des vari{\'e}t{\'e}s et ce bien que la structure de ce cone soit notamment importante dans le cadre des conjectures de Manin sur le d{\'e}compte de points rationnels (voir par exemple \cite{Peyre}).

\noindent  La r{\'e}solution canonique nous permet de relier la structure du cone effectif de la surface construite aux valeurs de l'indice effectifs:

 La valeur de l'indice effectif ainsi que la connaissance du diviseur canonique de la surface $V$ o{\`u} $(\pi,V)$ est la r{\'e}solution canonique nous permettent 
 de donner des informations sur la structure du cone effectif de la surface associ{\'e}e {\`a} la r{\'e}solution canonique.
Avec les notations du paragraphe 4, le diviseur canonique de $V$ s'exprime sous la forme:
 il  existe $a_i,a'_j\in\mathbb{Z}$ pour $1\leq i\leq n$ et $1\leq j\leq m$ tels que :
\vspace{-0.8cm} 

$$K_V=-3H^{\sharp}+\sum_{i=1}^{n}a_i E_i+ \sum_{j=1}^{m}a'_jF_j.$$

\paragraph*{{\bf{Th{\'e}or{\`e}me  B}}}
{\it{Soit $\phi\in \Aut(\mathbb{A}^2)$ de degr{\'e} alg{\'e}brique au moins deux.
 Alors, pour la r{\'e}solution canonique $(\pi,V)$, nous avons:
\begin{itemize}
\item si $\alpha(\phi, eff)=\frac{-1}{3}(a_n+a'_m)$ le cone des diviseurs effectifs de $V$ est poly{\'e}dral, 
\item si $\alpha(\phi, eff)>\frac{-1}{3}(a_n+a'_m)$ le cone des diviseurs effectifs de $V$ n'est pas poly{\'e}dral.
\end{itemize}}}

\noindent Le plan de ce texte est le suivant, nous 
faisons des rappels sur le th\'eor\`eme g\'en\'eral de structure du cone effectif d'une surface, ainsi que sur les automorphismes de l'espace affine. Puis nous construisons la r\'esolution canonique associ\'ee \`a un automorphisme du plan affine, nous \'etudions ces propri\'et\'es g\'eom\'etriques et enfin nous prouvons les th\'eor\`eme A et B.

\paragraph*{Remerciements}  Ce travail a \'et\'e effectu\'e au Max-Planck Institut fuer Mathematik que je remercie pour les excellentes conditions de travail.

\section{Structure du cone effectif d'une surface}
La structure du cone effectif d'une surface lisse projective d{\'e}coule notamment du th{\'e}or{\`e}me de l'indice de Hodge.
En effet, d'apr{\`e}s le th{\'e}or{\`e}me de l'indice de Hodge, il existe une base de NS(S) pour laquelle la forme d'intersection est donn{\'e}e par :
$$x_1^2-\sum_{i=2}^m x_i^2.$$
Consid{\'e}rons,
$$Q^+:=\{x\in NS(S)\ \mbox{tels que}\ x_1>(\sum_{i=2}^m x_i^2)^{\frac{1}{2}}\}.$$
La structure du cone effectif est alors la suivante.

\begin{theo}[\cite{Kollar}II.4.13]\label{conedim2}
Soit S une surface projective lisse.
Alors, 
$$\Pic^+_{\mathbb{R}}(S)=\overline{\mbox{Q}}^{+}+\sum_{\mbox{D}} \mathbb{R}^{+}[\mbox{D}],$$
 o{\`u} la sommation se fait sur toutes les courbes irr{\'e}ductibles 
$\mbox{D}\subset S$ telles que  $\mbox{D}^2<0.$
\end{theo}
Ainsi les ar{\^e}tes extr{\'e}males sont les courbes irr{\'e}ductibles d'auto-intersection n{\'e}gative.

\begin{rema} \label{KS} De mani{\`e}re {\'e}quivalente on peut voir $\overline{\mbox{Q}}^{+}$ comme  $ \Pic^+_{\mathbb{R}}(S)_{K_{S}\geq 0}$ o{\`u} $K_S$ d{\'e}signe le diviseur canonique de $S$.
\end{rema}

\section{Propri{\'e}t{\'e}s et d{\'e}finition des automorphismes affines}

\begin{lemm}\label{rat}(\cite{Beauville} lemme II.10)
Soit $\phi:S\cdots\rightarrow S'$ une application rationnelle entre deux surfaces avec $\phi^{-1}$ non d{\'e}finie en un point $P\in S'$.
Il existe une courbe $C$ de $S$ telle que:
$$\phi(C)=P.$$
\end{lemm}

\begin{lemm}\label{aff1}(\cite{Sibony} p.124)
Soit $\phi$ un automorphisme non-lin{\'e}aire de $\mathbb{A}^r$ . Nous avons alors\,:\\
$$\phi(H\setminus Z(\phi))\subset Z(\phi^{-1})\qquad
{\mbox {et}} \qquad 
\phi^{-1}(H\setminus Z(\phi^{-1}))\subset Z(\phi).$$
\end{lemm}

La notion suivant a \'et\'e introduite par N. Sibony \cite{Sibony}.

\begin{defi}
Soit $\phi$ un automorphisme de $\mathbb{A}^r$ avec $r \geq 2$. L'automorphisme $\phi$
est dit {\it{r{\'e}gulier}} 
si $\deg(\phi)>1$ et \,:
$$Z(\phi)\cap Z(\phi^{-1})=\emptyset.$$
\end{defi}

La d{\'e}finition suivante reprend le vocabulaire d{\'e}fini par S. Lamy \cite{lamy}.

\begin{defi}\label{propre}
Soit $\phi\in \Aut(\mathbb{A}^2)$. Soit $\pi:V\rightarrow \mathbb{P}^2$ un produit fini d'{\'e}clatements tel que $\phi\pi$ soit un morphisme. Nous appellerons {\it{points d'ind{\'e}termination}} de $\phi$ les points que l'on {\'e}clate lors de la construction de  $\pi$. Ces points appartiennent {\`a}  $\mathbb{P}^2$ ou {\`a} des {\'e}clat{\'e}s de  $\mathbb{P}^2$.
Les points d'ind{\'e}termination situ{\'e}s dans $\mathbb{P}^2$ seront dits {\it{points d'ind{\'e}termination propres}}.
\end{defi}

\begin{lemm} (\cite{lamy} lemme 9)\label{lamylem}
Soit $\phi\in \Aut(\mathbb{A}^2)$ de degr{\'e} au moins 2. Alors:
\begin{enumerate}
\item $\phi$ admet un seul point d'ind{\'e}termination propre situ{\'e} sur $H$.

\item $\phi$ admet des  points d'ind{\'e}termination $P_1,\dots,P_s$ ($s\geq 2$) tels que :
\begin{enumerate}
\item  $P_1$ soit le point d'ind{\'e}termination propre,
\item pour tout $i=2,\dots,s$ le point $P_i$ soit situ{\'e} sur le diviseur produit en {\'e}clatant $P_{i-1}$.
\end{enumerate}
\end{enumerate}
\end{lemm}

Le lemme suivant \ref{lamy1} d{\'e}coule imm{\'e}diatement de la premi{\`e}re assertion du lemme \ref{lamylem}.
\begin{lemm}\label{lamy1}
Soit $\phi\in \Aut(\mathbb{A}^2)$.
Alors,
$$\Card(Z(\phi)\cup Z(\phi^{-1}))\leq 2.$$
\end{lemm}
Plus pr{\'e}cis{\'e}ment dans ce m{\^e}me lemme 9 \cite{lamy}, nous avons, pour tout automorphisme non lin{\'e}aire $\phi$ du plan affine, il existe une suite finie d'{\'e}clatements $\pi_i:V_{i}\rightarrow V_{i-1}$ de centre $P_i$ avec $1\leq i\leq n$ avec $V_0=\mathbb{P}^2$, telle que:
$\phi\pi_1\cdots\pi_n$ est un morphisme de $V_n$ dans $\mathbb{P}^2$ et $\phi\pi_1\cdots\pi_{n-1}$ n'est pas un morphisme.

\noindent Les points $P_i$ avec $1\leq i\leq n$ sont des points infiniment voisins.

\section{{\'E}tude d'une r{\'e}solution canonique}
\subsection{Construction}
Nous utilisons ces propri{\'e}t{\'e}s g{\'e}om{\'e}triques des automorphismes du plan affine pour construire les r{\'e}solutions canoniques pour ces applications.
\subsubsection{Cas des automorphismes r{\'e}guliers}
Nous utiliserons leur d{\'e}finition et le lemme \ref{lamy1}.

\noindent Soit $\phi\in \Aut(\mathbb{A}^2)$ un automorphisme r{\'e}gulier. Nous avons deux familles distinctes de points associ{\'e}es {\`a} $\phi$ et $\phi^{-1}$, il s'agit de :
$$P_1,\dots,P_n$$  
$$\mbox{et}\quad Q_1,\dots,Q_m$$ 
associ{\'e}e respectivement {\`a} $\phi$ et {\`a} $\phi^{-1}$. Quitte {\'e}changer $\phi$ et $\phi^{-1}$, nous pouvons supposer $n\leq m$.
Nous d{\'e}finissons $\pi_i:V_{i}\rightarrow V_{i-1}$ l'{\'e}clatement de centre $P_i$ et $Q_i$ pour $1\leq i\leq n$ et de centre $Q_i$ pour $n\leq i\leq m$. La vari{\'e}t{\'e} $V_m$ sera not{\'e}e $V$, et $\pi_1\cdots\pi_m$ sera not{\'e} $\pi$.
Les diviseurs exceptionnels associ{\'e}s seront not{\'e}s:
$$E_1,\dots ,E_n$$ 
$$F_1,\dots, F_m.$$
La r{\'e}solution est minimale et unique, les entiers $n$ et $m$ {\'e}tant choisis comme dans le paragraphe pr{\'e}c{\'e}dent.
\subsubsection{Cas des automorphismes non-r{\'e}guliers}Nous utiliserons le lemme \ref{lamy1}.

\noindent Soit $\phi\in \Aut(\mathbb{A}^2)$ un automorphisme non-r{\'e}gulier non-lin{\'e}aire.

\noindent Il existe $i_0 \in \mathbb{N}\setminus\{0\}$ tel que $\pi_i:V_{i}\rightarrow V_{i-1}$ pour $1\leq i\leq i_0$ des {\'e}clatements de centre respectivement $P_1,\dots,P_{i_0}$, nous noterons
$$\psi_i=\phi\pi_1\cdots\pi_{i}\quad  \mbox{et} \quad \psi'_i=\phi^{-1}\pi_1\cdots\pi_{i},$$ 

\noindent avec pour tout $1\leq i\leq i_0$, $Z(\psi_i)=Z(\psi'_i)$ et $ P_{i_0+1}=Z(\psi_{i_0})\neq Z(\psi'_{i_0})=Q_{i_0+1}$.

\noindent Nous avons donc deux familles de points infiniment voisins: $P_1,\dots,P_{i_0},\dots, P_n$ et  $P_1,\dots,P_{i_0},Q_{i_0+1}\dots,Q_m$. 

\noindent Nous pouvons supposer $n\leq m$. Les {\'e}clatements pour $i\geq i_0+1$ sont  d{\'e}finis de la m{\^e}me mani{\`e}re que pour le cas r{\'e}gulier.
Les diviseurs exceptionnels associ{\'e}s sont:
$$E_1,\dots,E_n$$ 
$$E_1,\dots ,E_{i_0},F_{i_0+1},\dots, F_m.$$
Afin de traiter, dans la mesure du possible, simultan{\'e}ment les deux cas nous noterons la seconde famille de diviseurs exceptionnels :
$$F_1,\dots ,F_{i_0},F_{i_0+1},\dots, F_m,$$
avec pour tout $1\leq i\leq i_0$, $E_i=F_i$.

\subsection{Premiers r{\'e}sultats}

\paragraph{Notations} Soit $\phi \in\Aut(\mathbb{A}^2)$ non lin{\'e}aire. Soit $(\pi,V)$ une r{\'e}solution canonique de $\phi$. Soit $H^{\sharp}$ la transform{\'e}e stricte de $H$ et $E_1,\dots E_n,F_1,\dots,F_m$ les diviseurs exceptionnels associ{\'e}s {\`a} cette r{\'e}solution.  
Le groupe de Picard $\Pic(V)$ est donc engendr{\'e} par $H^{\sharp},E_1,\dots ,E_n,F_1,\dots, F_m$. Nous noterons $\psi=\phi\pi$ et $\psi=\phi^{-1}\pi$  les morphismes obtenus.

\noindent Il existe donc $b,c,e,b_i,c_i,e_i,b'_j,c'_j,e'_j\in \mathbb{Z}$ pour $1\leq i\leq n$ et $1\leq j\leq m$ tels que :
$$\psi^*(H)=bH^{\sharp}+\sum_{i=1}^{n}b_iE_i+\sum_{j=1}^{m}b'_jF_j,$$
$$\psi'^*(H)=cH^{\sharp}+\sum_{i=1}^{n}c_iE_i+\sum_{j=1}^{m}c'_jF_j,$$
$$\pi^*(H)=eH^{\sharp}+\sum_{i=1}^{n}e_iE_i+\sum_{j=1}^{m}e'_jF_j.$$
Si $\phi$ n'est pas r{\'e}gulier, alors il existe $i_0$ tel que pour tout $1\leq j\leq i_0$ nous avons :$b'_j=c'_j=e'_j=0$.

\begin{rema}
Dans tous les exemples que nous avons {\'e}tudi{\'e}s, nous avons  $n=m$.
\end{rema}

\begin{prop}\label{directe}
Soit $\phi \in\Aut(\mathbb{A}^2)$ non-lin{\'e}aire. Soit $(\pi,V)$ la r{\'e}solution canonique de $\phi$.
Les assertions suivantes sont  v{\'e}rifi{\'e}es:
\begin{enumerate}
\item Pour toute courbe irr{\'e}ductible $C\notin \{H^{\sharp},E_1,\dots E_n,F_1,\dots,F_m\}$, nous avons $\dim({\overline{\psi(C)}})=\dim({\overline{\psi'(C)}})=1$, ${\overline{\psi(C)}}\neq H$ et ${\overline{\psi'(C)}}\neq H$.
\item Nous avons $\psi_*(H^{\sharp})=\psi'_*(H^{\sharp})=0$.
\item Pour tout $1\leq i\leq n-1$ et pour tout $1\leq j\leq m$, nous avons :
$$\psi_*(E_{i})=0\quad \mbox{et}\quad \psi_*(F_{j})=0.$$
Pour tout $1\leq i\leq n$ et pour tout $1\leq j\leq m-1$, nous avons :
$$\psi'_*(E_{i})=0\quad \mbox{et}\quad \psi'_*(F_{j})=0,$$
de plus,
$$\psi_*(E_{n})=H\quad \mbox{et}\quad \psi'_*(F_{m})=H.$$
\end{enumerate}
\end{prop}

\begin{proof}
Nous consid\'erons chacun des points s\'epar\'ement.
\begin{enumerate}
\item Soit $C\notin \{H^{\sharp},E_1,\dots E_n,F_1,\dots,F_m\}$. Nous avons:
 $$\pi(C\setminus((\cup_{i=1}^{n}C\cap E_i)\cup (\cup_{j=1}^{m}C\cap F_j)))\not\subset H\setminus(P_0\cup Q_0),$$ 
d'o{\`u} le r{\'e}sultat car $\phi$ est construit {\`a} partir d'un automorphime de l'espace affine.
\item Nous avons $\pi(H^{\sharp}\setminus(\cup_{i=1}^{n}H^{\sharp}\cap E_i))= H\setminus(Z(\phi))$.

\noindent De l{\`a} 

$$\psi(H^{\sharp}\setminus(\cup_{i=1}^{n}H^{\sharp}\cap E_i))=\phi(H\setminus(Z(\phi))\subset Z(\phi^{-1}),$$ 
d'apr{\`e}s le lemme \ref{aff1}. Or, $\dim(Z(\phi^{-1}))=0$, donc $\dim({\overline{\psi(H^{\sharp})}})<1$, d'o{\`u} $\psi_*(H^{\sharp})=0$.
Le raisonnement est analogue pour $\psi'_*(H^{\sharp})$.
\item Nous noterons de la m{\^e}me mani{\`e}re $E_{n-1}$ et $\pi_n(E_{n-1})$. 
  Nous avons $\psi=\psi_{n-1}\pi_n$, or $\psi_{n-1}$ est une application rationnelle, donc d'apr{\`e}s le lemme \ref{rat} il existe $C$ une courbe de $\mathbb{P}^2$ telle que $\psi_{n-1}^{-1}(C)=P_n$, or par construction cette courbe ne peut {\^e}tre que $H$.

\noindent De plus, par construction $\psi_{n-1}(E_{n-1})\subset H$. Supposons que  ${\overline{\psi_{n-1}(E_{n-1})}}=H$, alors:

\centerline{Il existe un point $P\in E_{n-1}$ avec $P\neq P_0$ tel que $\psi_{n-1}(P)=Q\in Z(\psi_{n-1}^{-1});$}

de l{\`a} $\psi_{n-1}^{-1}(Q)=P=P_n$, d'o{\`u} la contradiction. Pour $1\leq i <n-1$ le raisonnement est analogue.

\noindent Pour $1\leq j\leq m$, si $\phi$ est r{\'e}gulier nous avons $\pi(F_j)=Q_1\notin Z(\phi)$ d'o{\`u} le r{\'e}sultat, si $\phi$ n'est pas r\'egulier le raisonnement est essentiellement le m{\^e}me si ce n'est que l'on consid{\`e}re $\psi_{i_1}$ qui est r{\'e}gulier. 

\noindent Enfin, comme $\psi$  est une application birationnelle il existe $C$ une courbe de $V$ telle que  ${\overline{\psi(C)}}=H$ et la seule courbe possible est $E_n$.

\end{enumerate}
\end{proof}

\begin{lemm}\label{relat}
Soit $\phi\in \Aut(\mathbb{A}^2)$ de degr{\'e} alg{\'e}brique au moins deux. Soit $(\pi,V)$ la r{\'e}solution canonique de $\phi$.
Avec les notations pr{\'e}c{\'e}dentes nous avons:
$$\psi^*(H).\psi'^*(H)=c_n=b'_m,$$
$$\psi^*(H).\pi^*(H)=e_n=b,$$
$$\psi'^*(H).\pi^*(H)=e'_m=c.$$
\end{lemm}

\begin{proof}
D'apr\`es la formule de projection, nous avons:
$$\psi^*(H).\psi'^*(H)=H.\psi_*\psi'^*(H)=c_n H^2=c_n,$$
et nous avons {\'e}galement:
 $$\psi^*(H).\psi'^*(H)=\psi'_*\psi^*(H).H=b'_m H^2=b'_m.$$
Les autres \'egalit\'es s'obtiennent de fa{\c{c}}on analogue.
\end{proof}

\paragraph{Diviseur canonique associ{\'e} {\`a} la r{\'e}solution canonique}Soit $(\pi,V)$ la r{\'e}solution canonique de $\phi$. Soit $K_V$ le diviseur canonique de $V$. 
\noindent Nous utilisons les notations pr{\'e}c{\'e}dentes. Par construction et d'apr{\`e}s \cite{hart} chapitre 5 proposition 3.3, il  existe $a_i,a'_j\in\mathbb{Z}$ pour $1\leq i\leq n$ et $1\leq j\leq m$ tels que :
$$K_V=-3H^{\sharp}+\sum_{i=1}^{n}a_i E_i+ \sum_{j=1}^{m}a'_jF_j;$$ 
en effet $-3H$ est le diviseur canonique de $\mathbb{P}^2$.

\begin{rema}
Si $\pi=\pi_1$, alors la surface $V$ obtenue est une surface de del Pezzo, et par cons{\'e}quent d'apr{\`e}s 
\cite{Kollar} chapitre 2 exemple 4.15.2, le cone des diviseurs effectifs est poly{\'e}dral de type fini.
\end{rema}

\subsection{Propri{\'e}t{\'e}s du diviseur $D(\alpha, \pi, V)$}

\begin{lemm}\label{propd}
Soit $\phi\in \Aut(\mathbb{A}^2)$ de degr{\'e} alg{\'e}brique au moins deux. Soit $(\pi,V)$ une r{\'e}solution canonique de $\phi$. Pour tout  $\alpha\in \mathbb{R}$, nous avons:
\begin{itemize}
\item le nombre d'intersection suivant:
$$D(\alpha,\pi,V).H^{\sharp}=-\alpha,$$
\item Pour tout $1\leq i\leq n-1$ et $1\leq j\leq m-1$ :
$$D(\alpha,\pi,V).E_i=D(\alpha,\pi,V).F_j=0,$$
$$D(\alpha,\pi,V).E_n=D(\alpha,\pi,V).F_m=1,$$
\item l'auto-intersection suivante:
$$D(\alpha,\pi,V)^2=2(1+c_n)-2\alpha(b+c)+\alpha^2.$$
\end{itemize}
\end{lemm}

\begin{proof}
Nous avons :
\begin{eqnarray*}
D(\alpha,\pi,V)^2&=& \psi^*(H)^2 +\psi^*(H).\psi'^*(H)-\alpha \psi^*(H).\pi^*(H) \\
& & +\psi'^*(H).\psi^*(H)+\psi'^*(H)^2-\alpha \psi'^*(H).\pi^*(H)\\
& & -\alpha\pi^*(H).\psi^*(H)-\alpha\pi^*(H).\psi'^*(H)+\alpha^2\pi^*(H)^2\\
&=& H^2+2c_n-2\alpha b+ H^2+-2\alpha c +\alpha^2 H^2,\\ 
\end{eqnarray*}
d'o{\`u} le r{\'e}sultat.
\end{proof}

\section{Preuve des th\'eor\`eme A et B}

\subsection{R\'esultats pr\'eliminaires}

\begin{theo}[Th{\'e}or{\`e}me de Kleiman (\cite{Lazarsfeld} th 1.4.8.)]\label{klei}
Soit $X$ une vari{\'e}t{\'e} (un sch{\'e}ma) complet. Si $D$ est un $\mathbb{R}$-diviseur nef de $X$, alors pour toute sous-vari{\'e}t{\'e} (sch{\'e}ma) irr{\'e}ductible $V\subset X$ de dimension $K$, nous avons :
$$D^K.V\geq 0.$$ 
\end{theo}

\begin{theo}(\cite{Dema} 6.6) \label{dema}
Soit $X$ une vari\'et\'e alg\'ebrique projective. Nous avons,
$$ \overline  \Pic^{a}(X)=\Pic^{nef}(X).$$
\end{theo}

\subsection{Calculs des nombres d'intersection}
\begin{prop}\label{dapv}
Soit $\phi\in \Aut(\mathbb{A}^2)$ de degr{\'e} alg{\'e}brique au moins deux. Soit $(\pi,V)$ la r{\'e}solution canonique de $\phi$. Soit $\alpha\leq \alpha(\phi,eff)$ alors, nous avons :
\begin{itemize}
\item si $\alpha<\frac{-1}{3}(a_n+a'_m)$ alors  $D(\alpha,\pi,V)\in \Poly(V),$
 o{\`u} $\Poly(V)$ d{\'e}signe la partie poly{\`e}drale du cone effectif de $V$,
\item si $\alpha=\frac{-1}{3}(a_n+a'_m) $ alors  $D(\alpha,\pi,V)\in \Pic^+_{K_{V}=0}(V)$,
\item si $ \frac{-1}{3}(a_n+a'_m)<\alpha$  alors  $D(\alpha,\pi,V)\in \Pic^+_{K_{V}>0}(V)$.
\end{itemize}
\end{prop}

\begin{proof}
Soit $(\pi,V)$ la r{\'e}solution canonique de $\phi$. Soit $K_V$ le diviseur canonique de $V$.

Nous calculons le nombre d'intersection de $D(\alpha,\pi,V)$ avec le diviseur canonique pour d{\'e}terminer dans quelle partie du cone effectif se trouve $D(\alpha,\pi,V)$.
{\`A} l'aide du lemme \ref{propd}, nous obtenons.
\begin{eqnarray*}
D(\alpha,\pi,V).K_V&=& \psi^*(H)+\psi'^*(H)-\alpha  \pi^*(H)).K_V\\
&=& \psi^*(H).K_V+\psi'^*(H).K_V-\alpha  \pi^*(H).K_V\\
&=&H.\psi_*(K_V)+H.\psi'_*(K_V)-\alpha  H.\pi_*(K_V) \\
&=&H.(a_nH)+H.(a'_mH)-\alpha  H.(-3H)\\
&=&a_n+a'_m+3\alpha\\
\end{eqnarray*}
D'o{\`u} le r{\'e}sultat d'apr{\`e}s le th{\'e}or{\`e}me \ref{conedim2}.
\end{proof}

La preuve du th\'eor\`eme A d\'ecoule imm\'ediatement des th\'eor\`emes  \ref{klei} et \label{dema} et du premier point du lemme \ref{propd}.

La preuve du th{\'e}or{\`e}me B repose essentiellement sur la proposition \ref{dapv} et la d{\'e}finition de l'indice effectif.

\begin{rema}
Soit $\phi\in \Aut(\mathbb{A}^2)$. S'il existe une r{\'e}solution $(\pi,V)$ de $\phi$ pour laquelle $V$ est une surface de del Pezzo 
alors nous avons $\alpha(\phi,eff)\leq \frac{4}{3}$.
Ce cas se produit notamment, si suivant la terminologie utilis{\'e}e dans \cite{lamy}, $\phi$ et $\phi^{-1}$ ont un point d'ind{\'e}termination propre (voir d\'efinition \ref{propre}) et ce que ces points soient confondus ou non. 
\end{rema}

\section{Valeurs num{\'e}riques de l'indice effectif}

Dans ce paragraphe, nous donnons les valeurs des deux invariants  g{\'e}om{\'e}triques pour diff{\'e}rents exemples, dans tous ces exemples $a$ est une constante non nulle. J. Silverman a d{\'e}termin{\'e} la valeur de l'indice effectif pour l'un de ces exemples, l'application de H{\'e}non  g{\'e}n{\'e}ralis{\'e}e (\cite{Sil1}). Pour calculer ces indices nous avons utilis{\'e} la r{\'e}solution canonique. 
Le d{\'e}tail des calculs est omis.

Nous rappelons la d\'efinition du degr\'e dynamique (voir par exemple \cite{Sibony}): 

\begin{defi}
Soit $\phi$ un automorphisme affine. La limite suivante existe et
d{\'e}fini le degr{\'e} dynamique $\delta(\phi)$:
  $$\delta(\phi):=\inf_{n\geq 1}
deg(\phi^n)^{1/n}.$$
\end{defi}

%\subsection{D{\'e}termination de l'indice ample}
%Nous utilisons les calculs effectu{\'e}s pour la d{\'e}termination de l'indice effectif.

%\noindent Soit $B(H^{\sharp},E_i)$ avec $i\in I\subset \mathbb{N}\setminus \{0\}$ et $I$ ensemble fini non vide. Par construction, $B$ est une base de $\Pic(V)$.
%Pour utiliser le crit{\`e}re de Nakai-Moishezon (\cite{hart} chap.5 th. 1.10), nous devons consid{\'e}rer toutes les courbes irr{\'e}ductibles de $V$, dans notre cas il s'agit de v{\'e}rifier ce crit{\`e}re avec les {\'e}l{\'e}ments de $B$.

%\noindent Dans chacun des exemples consid{\'e}r{\'e}s, $\exists i_0\in I$ tel que 
%$$D(\alpha(\phi, eff),\pi,V).E_{i_{0}}=0.$$
%Le diviseur $D(\alpha(\phi, eff),\pi,V)$ n'est par cons{\'e}quent pas un diviseur ample, ce qui est  {\'e}galement le cas pour tous les el{\'e}ments de sa classe d'{\'e}quivalence num{\'e}rique, ceci d{\'e}coule directement du th{\'e}or{\`e}me de Nakai-Moishezon.

%\noindent Par contre, gr{\^a}ce {\`a} ces m{\^e}me calculs, nous obtenons que  $D(\alpha(\phi, eff),\pi,V)$ est un diviseur nef (num{\'e}riquement effectif). Or $\Pic^{a}(V)\subset \Pic^{nef}(V)$ o{\`u} $\Pic^{nef}$ d{\'e}signe le cone des diviseurs nef, de plus le cone ample est l'int{\'e}rieur du cone nef ( voir \cite{Lazarsfeld} th. 1.4.21 ou \cite{Dema} prop 6.6, d'o{\`u} les r{\'e}sultats que nous obtenons.

\vspace{1cm}

\noindent\begin{tabular}{|c||c|c|c|c|c|}
\hline
Application & degr{\'e} & degr{\'e}  & indice  \\
 &  alg{\'e}brique &  dynamique & effectif \\
\hline
$\phi(x,y)=$& $\deg(\phi)$ & $\delta(\phi)$ & $\alpha(\phi,eff)$\\
\hline
\hline
%$(x+ay^2+by+c,y)$&2  & 1 & $\frac{3}{2}=1+\frac{1}{2}$& $\frac{3}{2}$& $\leq 0$\\
%\hline
%$(x+ay^3,y)$&3 & 1 & $\frac{4}{3}=1+\frac{1}{3}$& $\frac{4}{3}$& $\leq 0$ \\
%\hline
%\hline
$(y,y^2+b+ax)$&2   & 2& $ \frac{5}{2}=2+\frac{1}{2}$ \cite{Sil1}\\
\hline
$(y+ax^3,x)$ &3 & 3 & $\frac{10}{3}=3+\frac{1}{3}$\\
\hline
$(y+ax^4,x)$&4 &4& $\frac{17}{4}=4+\frac{1}{4}$\\
\hline 
\end{tabular}

\vspace{1cm}

\begin{rema} Pour les exemples de degr{\'e} dynamique diff{\'e}rents de 1, nous obtenons la borne maximale du th{\'e}or{\`e}me C.
Sur les calculs effectu{\'e}s, nous observons :
$\alpha(\phi,eff)=\delta(\phi)+\frac{1}{\deg(\phi)}.$
Une question naturelle se pose, les indices amples et effectifs ne seraient-il pas rationnels?
\end{rema}


\begin{thebibliography}{9}








\bibitem[B83]{Beauville}
A. Beauville, Complex algebraic surfaces, London Mathematical Society Lecture Note series 68 (1983).


\bibitem[Dem96]{Dema}
Jean-Pierre Demailly, $L\sp 2$ vanishing theorems for positive line bundles and adjunction theory. Transcendental methods inalgebraic geometry (Cetraro, 1994), 1--97, Lecture Notes in Math., 1646, Springer, Berlin, (1996). 


\bibitem[Den95]{Denis}
L. Denis, Points p{\'e}riodiques des automorphismes affines, J. reine angew. Math. 467, 157-167 (1995).



\bibitem[H77]{hart}
R. Hartshorne, Algebraic geometry, Graduate Texts in Mathematics  52, Springer Verlag (1977).

\bibitem[HS00]{hs}
M. Hindry, J. Silverman, Diophantine geometry, an introduction, Graduate Texts in Mathematics 201, Springer Verlag (2000).


\bibitem[K96]{Kollar}
J. Koll{\'a}r, Rational curves on algebraic varieties, Ergebnisse des Mathematik und ihrer Grenzgebiete 3. Folge.Band 32,
 Springer (1996)


\bibitem[Lam02]{lamy}
S. Lamy, Une preuve g{\'e}om{\'e}trique du th{\'e}or{\`e}me de Jung, L'enseignement math{\'e}matique 48 291-315 (2002).

\bibitem[Laz01]{Lazarsfeld}
R. Lazarsfeld, Positivity in algebraic geometry, manuscrit (2001).


\bibitem[M1]{moi}
S. Marcello, Sur les propri{\'e}t{\'e}s arithm{\'e}tiques des it{\'e}r{\'e}s
d'automorphismes r{\'e}guliers, C. R. Acad. Sci. Paris, t.331, S{\'e}rie
I,11-16 (2000).

\bibitem[M2]{1art}
S. Marcello, Sur la dynamique arithm{\'e}tique des automorphismes de l'espace affine, Bulletin de la S.M.F.  31  229-257 (2003).

\bibitem[M3]{2art}
S. Marcello, Sur la dynamique arithm{\'e}tique $p$-adique des automorphismes de l'espace affine (2003).

\bibitem[M4]{moi4}
S. Marcello, G\'eom\'etrie, points rationnels et it\'er\'es des automorphismes de l'espace affine (2003).

\bibitem[P01]{Peyre}
E. Peyre, Points de hauteur born{\'e}e et g{\'e}om{\'e}trie des vari{\'e}t{\'e}s [d'apr{\`e}s Y. Manin et {\it{al.}}], S{\'e}minaire Bourbaki volume 2000-2001   Ast\'erisque 282  323-344  (2002).

\bibitem[S99]{Sibony}
N. Sibony, Dynamique des applications rationnelles de ${\mathbb{P}} ^k$, Panoramas et Synth{\`e}ses 8 (S.M.F.), 97-195 (1999)  .


\bibitem[Si94]{Sil1}
J. Silverman, Geometric and arithmetic properties of the H{\'e}non map, Math. Z. 215, 237-250 (1994).




\end{thebibliography}
\end{document}